\documentstyle[12pt,leqno,a4]{article}
\begin{document}
\footnotetext{AMS (1991) subject classification: Primary 39B12, 39B22.}
\large
$$ $$
\begin{center}
{\Large \bf On bounded solutions\\
of a problem of R. Schilling}
\end{center}
\vspace{1ex}
\begin{center}
{\bf Janusz Morawiec}
\end{center}
\vspace{1ex}
\begin{quotation}
\center
{\small
\begin{flushleft}
{\bf Abstract.} It is proved that if
\end{flushleft}
$0<q\leq (1-\sqrt[3]{2}+\sqrt[3]{4})/3,$
\vspace{-1.5ex}
\begin{flushleft}
then the zero function is the only solution
$f\! :\! {\rm I\! R}\rightarrow {\rm I\! R}$ of (1) satisfying (2) and
bounded in a neighbourhood of at least one point of the set (3).
\end{flushleft}}
\end{quotation}

The paper concerns bounded solutions
$f\! :\! {\rm I\! R}\rightarrow {\rm I\! R}$ of the functional equation
\begin{equation}
f(qx)=\frac{1}{4q}[f(x-1)+f(x+1)+2f(x)]
\end{equation}
such that
\begin{equation}
f(x)=0 \hspace{1cm} {\rm for} \hspace{.5cm} \mid x \mid >Q
\end{equation}
where $q$ is a fixed number from the open interval $(0,1)$ and
$$
Q=\frac{q}{1-q}.
$$

In what follows any solution $f\! :\!{\rm I\! R}\rightarrow {\rm I\! R}$
of (1) satisfying (2) will be called a {\it solution of Schilling's
problem}. In the present paper we are interested in bounded solutions of
Schilling's problem. The first theorem in this direction was obtained by
K.Baron in [1]. This theorem reads as follows:

{\it If $q\in(0,\sqrt{2}-1]$, then the zero function
is the only solution of Schilling's problem which is bounded in a
neighbourhood of the origin.}

This paper generalizes the above theorem in two directions.
\linebreak
Namely, the
interval $(0,\sqrt{2}-1]$ is replased by the larger one
\linebreak
$(0,\frac{1}{3}-\frac{\sqrt[3]{2}}{3}+\frac{\sqrt[3]{4}}{3}]$
and instead of the
boundedness in a neighbourhood of the origin we have boundedness in a
neighbourhood of at least one point of the set
\begin{equation}
\left\{\varepsilon\sum_{i=1}^{n}q^{i}\:\:\::\:\:n\in{\rm
I\!N}\cup\{0,+\infty\}\:,\:\varepsilon\in\{-1,1\}\right\}.
\end{equation}
$$ $$
(To simplify formulas we adopt the convention $\sum_{i=1}^{0}a_{i}=0$
for all real sequences $(a_{i}:i\in {\rm I\! N}).$) In other words, we
shall prove the following.

{\bf Theorem.}
{\it If
\begin{equation}
0<q\leq\frac{1}{3}-\frac{\sqrt[3]{2}}{3}+\frac{\sqrt[3]{4}}{3},
\end{equation}
then the zero function is the only solution of Schilling's problem which
is bounded in a neighbourhood of at least one point of the set} (3).

The proof of this theorem is based on two lemmas. However, we start
with the following simple remarks.

{\bf Remark 1.}
{\it If $f$ is a solution of Schilling's problem then so is the
function $g\! :\! {\rm I\! R}\rightarrow {\rm I\! R}$ defined by
the formula $g(x)=f(-x).$}

{\bf Remark 2.}
{\it Assume that $f$ is a solution of Schilling's problem.

{\rm (i)} If $q\neq\frac{1}{4}$, then $f(Q)=0$.
If $q=\frac{1}{4}$, then $f(Q)=0$ iff $f(qQ)=0.$

{\rm (ii)} If $q<\frac{1}{2}$, then $f(0)=0.$}

{\bf Proof.}
It is enough to put in (1): $x=Q/q$, $x=Q$ and $x=0$,
respectively, and to use condition (2).

{\bf Lemma 1.}
{\it Assume $q\in (0,\frac{1}{2}).$ If a solution of Schilling's problem
vanishes either on the interval $(-q,0)$ or on the interval $(0,q)$,
then it vanishes everywhere.}

{\bf Proof.}
Let $f$ be a solution of Schilling's problem vanishing on the
interval $[0,q)$. We shall prove
that $f$ vanishes on the interval $[0,Q)$. Define a sequence of sets
$(A_{n}:n\in{\rm I\!N})$ by the formula:
$$
A_{n}=[\;0\;,\;\sum_{i=1}^{n}q^{i}\;).
$$
Fix a positive integer $n$ and suppose that $f$ vanishes on the set
$A_{n}$. We shall show that $f$ vanishes also on the set $A_{n+1}$. To
this end fix an $x_{0}\in A_{n+1}\setminus A_{n}$. Putting
$x=x_{0}/q$ into (1) and taking into account that
$x-1\in A_{n}$, whereas $x+1>x>1>Q$ we get
\begin{equation}
f(x_{0})=\frac{1}{4q}[f(x-1)+f(x+1)+2f(x)]=0.
\end{equation}
Consequently $f$ vanishes on the set $\bigcup_{n=1}^{\infty }A_{n}$
which equals to $[0,Q)$. This and Remark 2 (i) show that $f$ vanishes on
$[0,+\infty).$ Hence and from (1) we infer that $f$ vanishes everywhere.

The case of the interval $(-q,0)$ reduces to the previous one by
Remark 1.

{\bf Lemma 2.}
{\it Assume $q\in (0,\frac{1}{2})$. If $f$ is a solution of Schilling's
problem, then
\begin{equation}
f(q^{m+n}x+\varepsilon\sum_{i=1}^{n}q^{i})=
\left(\frac{1}{2}\right)^{n}\left(\frac{1}{2q}\right)^{m+n}f(x)
\end{equation}
for all $x\in(Q-1,1-Q)$, for all $\varepsilon\in\{-1,1\}$, and for all
non-negative integers $m$ and $n$.}

{\bf Proof.}
Fix an $x_{0}\in(Q-1,1-Q)$. First we shall show that
\begin{equation}
f(q^{m}x_{0})=\left(\frac{1}{2q}\right)^{m}f(x_{0})
\end{equation}
for all non-negative integers $m$. Of course (7) holds for $m=0$.
Suppose that (7) holds for an $m$. Putting $x=q^{m}x_{0}$ into (1)
and using (2) and (7) we have
$$
f(q^{m+1}x_{0})=\frac{1}{4q}[f(x-1)+f(x+1)+2f(x)]=
\frac{1}{2q}f(x)
$$
$$
\hspace{-16ex}=\left(\frac{1}{2q}\right)^{m+1}f(x_{0}).
$$
This proves that (7) holds for all non-negative integers $m$.

Fix now a non-negative integer $n$ and suppose that (6) is satisfied
for all $x\in (Q-1,1-Q)$,for all $\varepsilon\in\{-1,1\}$, and
for all non-negative integers $m$.
Putting $x=q^{m+n}x_{0}+\varepsilon\sum_{i=1}^{n}q^{i}+\varepsilon$
into (1) and applying (2) and (6) with $x=x_{0}$ we obtain
$$
f(q^{m+n+1}x_{0}+\varepsilon\sum_{i=1}^{n+1}q^{i})=f(qx)=
\frac{1}{4q}[f(x-1)+f(x+1)+2f(x)]
$$
$$
\hspace{23ex}=\frac{1}{4q}f(x-\varepsilon)=
\left(\frac{1}{2}\right)^{n+1}\left(\frac{1}{2q}\right)^{m+n+1}f(x_{0}).
$$

The proof is completed.

Now we pass to the proof of the main theorem.

{\bf Proof of the theorem.}
It follows from (4) that $q<1/2$.

Fix $n\in {\rm I\!N}\cup\{0,+\infty\}$ and $\varepsilon\in\{-1,1\}$
such that a solution $f$ of Schilling's problem is bounded in a
neighbourhood of $\varepsilon\sum_{i=1}^{n}q^{i}.$ We may (and we do)
assume that $n$ is finite.

If $\mid x\mid <1-Q$ is fixed, then the left-hand-side of (6) is
bounded with respect to $m$ whereas $\lim_{m\rightarrow
\infty}(1/2q)^{m+n}=+\infty.$ This shows that
\begin{equation}
f(x)=0  \hspace{1cm} {\rm for} \hspace{.5cm} \mid x \mid<1-Q.
\end{equation}

Consider two cases:\\
(i) \hspace{2cm} $q\leq\frac{3-\sqrt{5}}{2}$\\
and\\
(ii) \hspace{1.9cm} $\frac{3-\sqrt{5}}{2}<q\leq
\frac{1}{3}-\frac{\sqrt[3]{2}}{3}+\frac{\sqrt[3]{4}}{3}.$

In the case (i) we have $q\leq 1-Q$ which jointly with (8) and
Lemma~1 gives $f=0$.

So we assume now that (ii) holds. First we notice that putting $x=1-Q$
into (1) and applying (8), Remarks 1 and 2 (i) and (2) we get
$$
0=f(q(1-Q))=\frac{1}{4q}[f(-Q)+f(2-Q)+2f(1-Q)]=\frac{1}{2q}f(1-Q).
$$
Hence, from (8) and Remark~1 we obtain
\begin{equation}
f(x)=0  \hspace{1cm} {\rm for} \hspace{.5cm} \mid x \mid\leq 1-Q.
\end{equation}

Fix an $x_{0}\in[qQ,q(2-Q)].$
Putting $x=x_{0}/q$ into (1) and using (9), (2) and Remark~2
we have (5). Similarly (cf. Remark 1), $f(x)=0$ for $x\in[-q(2-Q),-qQ]$.
Consequently,
\begin{equation}
f(x)=0\hspace{1cm}{\rm whenever}\hspace{.5cm}
qQ\leq\mid x\mid\leq q(2-Q).
\end{equation}

Now we fix an $x_{0}\in[q-q^{2}(2-Q),q^{2}(2-Q)].$ Putting
$x=x_{0}/q$ into (1), taking into account the inequality
$qQ<1-q(2-Q)$ and applying (10) and (2) we obtain (5) once again.
Similarly $f(x)=0$ for $x\in[-q^{2}(2-Q),-q+q^{2}(2-Q)]$ and so
\begin{equation}
f(x)=0  \hspace{1cm} {\rm whenever} \hspace{.5cm}
q-q^{2}(2-Q)\leq\mid x \mid\leq q^{2}(2-Q).
\end{equation}
As the function $3t^{3}-3t^{2}+3t-1$ increases and vanishes at
\linebreak
$(1-\sqrt[3]{2}+\sqrt[3]{4})/3,$ we have
\begin{equation}
q-q^{2}(2-Q)\leq 1-Q.
\end{equation}
Relations (9), (12) and (11) give
\begin{equation}
f(x)=0 \hspace{1cm}{\rm for}\hspace{.5cm}\mid x \mid\leq q^{2}(2-Q).
\end{equation}

Now let us fix an $x_{0}\in[1-q(2-Q),1-qQ].$ Putting $x=x_{0}-1$ into
(1) and using (13), (2) and (10) we have
\begin{equation}
0=f(qx)=\frac{1}{4q}[f(x-1)+f(x+1)+2f(x)]=\frac{1}{4q}f(x_{0}).
\end{equation}
So we obtain
\begin{equation}
f(x)=0 \hspace{1cm}{\rm whenever}\hspace{.5cm}
1-q(2-Q)\leq x \leq 1-qQ.
\end{equation}
Since (cf. (12)) $q+q^{2}Q<1-qQ$ and $1-q(2-Q)<q(2-Q)$, (15) proves that
\begin{equation}
f(x)=0 \hspace{1cm}{\rm whenever}\hspace{.5cm}
q(2-Q)\leq x \leq q+q^{2}Q.
\end{equation}

Finally assume that  $1-Q\leq x_{0}\leq qQ.$ Putting $x=x_{0}+1$ into
(1) and using (16) and (2) we see that (14) holds. Hence
$$
f(x)=0 \hspace{1cm}{\rm whenever}\hspace{.5cm} 1-Q\leq x \leq qQ,
$$
which jointly with (9) and (10) gives
$$
f(x)=0 \hspace{1cm}{\rm whenever}\hspace{.5cm} 0\leq x\leq q(2-Q).
$$
In particular, since $q<q(2-Q)$, $f$ vanishes on the interval $(0,q).$
This jointly with Lemma~1 completes the proof.

\vspace{1cm}
{\bf Acknowledgement.} This research was supported by the Silesian
University Mathematics Department (Iterative Functional
\linebreak
Equations program).

\vspace{3ex}
\begin{flushleft}
{\bf Reference}
\end{flushleft}
[1] K. Baron, {\it On a problem of R.Schilling.} Berichte der Mathema-

tisch-statistischen Sektion in der Forschungsgesellschaft Joanne-

um--Graz, Bericht Nr. 286 (1988).

\vspace{3ex}
\begin{flushleft}
Instytut Matematyki,\\
Uniwersytet \'{S}l\c aski,\\
ul. Bankowa 14,\\
PL-40-007 Katowice
\end{flushleft}

\end{document}